\documentclass[12pt]{article}    
\usepackage{amsmath}
\usepackage{amssymb}
\usepackage{amsmath,amssymb,amsfonts,amsthm,latexsym}
\usepackage{booktabs}
\usepackage{array} 
\usepackage{lscape}
\usepackage{graphicx}
\usepackage{xy}
\usepackage{amsmath}
\usepackage{amssymb}
\input diagxy

\usepackage[colorlinks=true]{hyperref}
\usepackage{epsfig}
\usepackage{colordvi}
\usepackage{latexsym}
\usepackage{authblk}
\def\mysection{\setcounter{equation}{0}\section}

\newtheorem{thm}{Theorem}
 
\newtheorem{lem}{Lemma}
\newtheorem{prop}{Proposition}
\newtheorem{cor}{Corollary}
\newtheorem{defn}{Definition}

\newcommand{\beq}{\begin{equation}}
\newcommand{\eeq}{\end{equation}}
\newcommand{\beqr}{\begin{eqnarray}}
\newcommand{\eeqr}{\end{eqnarray}}

\begin{document}
	\title{\Large {{Quasi-uniform structures and functors} }}
    \author[$*,\;a$]{Minani Iragi}
    \author[$b$]{David Holgate}

	\affil[$^{a}$]{{\footnotesize Department of Mathematical Sciences, University of South Africa, P.O. Box 392, 003 Unisa, South Africa.}}
	
	\affil[$b$]{{\footnotesize Department of Mathematics and Applied Mathematics University of the Western Cape,   Bellville 7535, South Africa and Institute of Mathematics, Faculty of Mechanical Engineering, Brno  University of Technology Technická 2, 616 69 Brno, Czech Republic.}}
	
	\date{}
	\maketitle
	\def\thefootnote{\fnsymbol{footnote}}
	\setcounter{footnote}{0}
	
	\footnotetext{ 
	E-mail addresses: $^{a}$84miragi2016@mail.com, $^{b}$dholgate@uwc.ac.za\\
	$^{*}$Corresponding author.\\	
	The work of the first author is based on a research supported wholly by the  National Research Foundation of South Africa (Grant Numbers : 129519). The Second author also acknowledges the National Research Foundation of South Africa and the support from the European Social Fund and the state budget of the Czech Republic in the project no.CZ.02.2.69/0.0/0.0/16-027/0008371.}
 \begin{abstract}
 {\indent We study a number of categorical quasi-uniform structures induced by functors. We depart from a category $\mathcal{C}$ with a proper $(\mathcal{E}, \mathcal{M})$-factorization system, then define the continuity of a $\mathcal{C}$-morphism with respect to two syntopogenous structures (in particular with respect to two quasi-uniformities) on $\mathcal{C}$ and use it to describe the quasi-uniformities induced  by  pointed and copointed endofunctors of $\mathcal{C}$. In particular, we demonstrate that every quasi-uniformity on a reflective subcategory of $\mathcal{C}$ can be lifted to a coarsest quasi-uniformity on $\mathcal{C}$ for which every reflection morphism is continuous.}
 
 {\indent Thinking of categories supplied with quasi-uniformities as large ``spaces'', we generalize the continuity of $\mathcal{C}$-morphisms (with respect to a quasi-uniformity) to functors.  We prove that for an $\mathcal{M}$-fibration or a functor that has a right adjoint, we can obtain a concrete construction of the coarsest quasi-uniformity for which the functor is  $continuous$. The results proved are shown to yield those obtained for categorical closure operators. Various examples considered at the end of the paper illustrate our results.}
\end{abstract}

\\

\mysection{Introduction}
{\indent $ \; \; \; \;$  The introduction of categorical closure operators (\cite{dikranjan1987closure}) by Dikranjan and Giuli was the point of departure for study of topological structures on categories. This approach eventually motivated the introduction of categorical interior (\cite{MR1810290}) and neighbourhood (\cite{MR2838385}) operators. While the categorical interior operators were shown to be pleasantly related to neighbourhood operators, a nice relationship between closure and neighbourhood operators in a category was lacking until the categorical topogenous structures (\cite{MR0286814, iragi2016quasi}) were recently introduced. Indeed the conglomerate of categorical topogenous structures is order isomorphic to the conglomerate all neighbourhood operators and contains both the conglomerates of all interior and all closure operators as reflective subcategories.
	
{\indent  A natural generalization of the definition of a categorical topogenous structures leads to the concept of categorical syntopogenous structure which provides a convenient setting to investigate a quasi-uniform structure on a category. This is the point of departure in (\cite {MR4239728, iragi2019quasi}) where a categorical quasi-uniform structure is introduced and studied. Moreover, the use of syntopogenous structures allows the description of a quasi-uniformity as a family of categorical closure operators (see e.g \cite{holgate2019quasi}). A recent account of this relationship between quasi-uniformity and closure operators can be found in \cite{MR4239728, MR3641225} . 
		
{\indent The present paper aims to further study a categorical quasi-uniform structure. Considering a category $\mathcal{C}$ with a proper $(\mathcal{E}, \mathcal{M})$-factorization system, we show that for a syntopogenous structure $\mathcal{S}$ on $\mathcal{C}$ and an $\mathcal{E}$-pointed endofunctor $(F,\eta)$ of $\mathcal{C}$, there is a coarsest syntopogenous structure $\mathcal{S}^{F,\eta}$ on $\mathcal{C}$ for which every $\eta_{X}: X\longrightarrow FX$ is $(\mathcal{S}^{F,\eta}, \mathcal{S})$-continuous. Since a categorical quasi-uniformity is equivalent to a co-perfect syntopogenous structure and simple co-perfect syntopogenous structures are equivalent to idempotent closure operators (see e.g \cite{holgate2019quasi}), $\mathcal{S}^{F,\eta}$ allows us to construct the quasi-uniform structure and the closure operator induced by a pointed endofunctor.  In particular, we demonstrate that every quasi-uniformity $\mathcal{U}$ on a reflective subcategory of $\mathcal{C}$ can be lifted to a coarsest quasi-uniformity $\mathcal{U}^{F,\eta}$ on $\mathcal{C}$ for which every reflection morphism is $(\mathcal{U}^{F, \eta}, \mathcal{U})$-continuous. When applied to spaces, $\mathcal{U}^{F, \eta}$ turns out to describe initial structures induced by reflection maps. Dual results shall be obtained in the case of a copointed endofunctor.
			
{\indent For a functor $F : \mathcal{A}\longrightarrow \mathcal{C}$ and quasi-uniformities $\mathcal{U}$ and $\mathcal{V}$ on $\mathcal{A}$ and $\mathcal{C}$ respectively, we introduce the $(\mathcal{U}, \mathcal{V})$-continuity of $F$. It is shown that if $F$ is an $\mathcal{M}$-fibration or has a right adjoint, then one can concretely describe the coarsest quasi-uniformity $\mathcal{V}^{F}$ on $\mathcal{A}$ for which $F$ is $(\mathcal{V}^{F}, \mathcal{V})$-continuous. We then use the categorical co-perfect syntopogenous structures, to obtain a concrete description of the largest closure operator making $F$ continuous.}

{\indent In section 4, we describe categorical quasi-uniform structures induced by (co) pointed endofunctors, which we construct using the syntopogenous structures (Proposition $4.4$, Theorems $4.4$ and $4.9$). It is interesting to note that particular cases of these quasi-uniform structures correspond to the closure operators obtained by Dikranjan and Tholen in \cite{MR2122803} (chapter 5, Theorems $5.12$ and $5.12^{\ast}$). The study of continuity of functors with respect to two quasi-uniform structures and its use to describe the initial quasi-uniform structures induced by an $\mathcal{M}$-fibration
or a functor having a right adjoint (Proposition $5.4$, $5.7$ and $5.9$, Theorems $5.5$ and $5.8$) are devoted to section 5. Finally in section 6, we present a number of examples to illustrate the results obtained.


\section{ Preliminaries} 
{\indent  Our blanket reference for categorical concepts is \cite{MR2240597}. The basic facts on categorical closure operators used here can be found in \cite{MR2122803} or \cite{dikranjan1987closure}. For the categorical topogenous, quasi-uniform and syntopogenous structures, we use \cite{iragi2019quasi} and \cite{MR0286814}.
Throughout the paper, we consider a category $\mathcal{C}$ supplied with a proper $(\mathcal{E}, \mathcal{M})$-factorization system for morphisms. The category $\mathcal{C}$ is assumed to be 
$\mathcal{M}$-complete so that pullbacks of $\mathcal{M}$-morphisms along $\mathcal{C}$-morphisms and arbitrary $\mathcal{M}$-intersections of $\mathcal{M}$-morphisms exist and are again in $\mathcal{M}$. 
For any $X\in \mathcal{C}$, sub$X = \{m\in \mathcal{M}\;|\;\mbox{cod}(m) = X\}$. It is ordered as follows: $n\leq m$ if and only if there exists $j$ such that $m\circ j = n$. If $m\leq n$ and $n\leq m$ then they are isomorphic. We shall simply write $m = n$ in this case.  Sub$X$ is a (possibly large) complete lattice with greatest element $1_{X} : X\longrightarrow X$ and the least element $0_{X} : O_{X}\longrightarrow X$.}

{\indent Any $\mathcal{C}$-morphism, $f:X\longrightarrow Y$  induces an image/pre-image adjunction 
	$f(m)\leq n$ if and only if $m\leq f^{-1}(n)$ for all $n\in \;$sub$Y$, $m\in \;$sub$X$ with $f(m)$ the $\mathcal{M}$-component of the $(\mathcal{E}, \mathcal{M})$-factorization of $f\circ m$ while $f^{-1}(n)$ is the pullback of $n$ along $f$. 
	We have from  the image/pre-image adjunction that $f(f^{-1}(n))\leq n$ (with $f(f^{-1}(n)) = n$ if $f\in \mathcal{E}$ and $\mathcal{E}$ is pullback stable along $\mathcal{M}$-morphisms) and $m\leq f^{-1}(f(m))$  (with $m = f^{-1}(f(m))$ if $f\in \mathcal{M}$)
	for any $n\in \;$sub$Y$ and $m\in \;$sub$X$. }

Applying adjointness repeatedly we obtain the lemma below.

\lem\label{lem-latex-2e}\label{L22}
	Let 
	$$\bfig
	\square[X'`Y'`X`Y;f'`p'`p`f]
	
	\efig$$

	be a commutative diagram. Then for any subobject $n\in \;$sub$Y'$,  $p'(f'^{-1}(n))\leq f^{-1}(p(n))$.
\endlem

\begin{defn}\label{D154}
	$A\; pointed \; endofunctor$ of $\mathcal{C}$ is a pair $(F,\eta)$ consisting of a functor $F : \mathcal{C}\longrightarrow \mathcal{C}$ and a natural transformation $\eta : 1_{\mathcal{C}}\longrightarrow F$.
	\end{defn}
For any $\mathcal{C}$-morphism $f : X\longrightarrow Y$, $(F,\eta)$ induces the commutative diagram below. 
$$\bfig
\square[X`FX`Y`FY;\eta_{X}`f`Ff`\eta_{Y}]

\efig$$
	If each $\eta_{X}\in \mathcal{F}$ where $\mathcal{F}$ is a class of  $\mathcal{C}$-morphisms, then $(F,\eta)$ is $\mathcal{F}$-pointed. 
A copointed endofunctor of $\mathcal{C}$ is defined dually. 

\begin{defn}
	$A\;closure\;operator$ $c$ on $\mathcal{C}$ with respect to $\mathcal{M}$ is given by a
	family of maps \\$\{c_{X}$: sub$X\longrightarrow \;$sub$X\;|\; X\in \mathcal{C}\}$ such that:
	\begin{itemize}
		\item [$(C1)$]  $m\leq c_{X}(m)$ for all $m\in \;$sub$X$;
		\item [$(C2)$]  $m\leq n\Rightarrow c_{X}(m)\leq c_{X}(n)$ for all $m, n\in \;$sub$X$;
		\item [$(C3)$]   every morphism $f: X\longrightarrow Y$ is $c$-continuous, that is: $f(c_{X}(m))\leq c_{Y}(f(m))$ for all
		$m\in \;$sub$X$.
	\end{itemize}
\end{defn}
	We denote by CL$(\mathcal{C}, \mathcal{M})$ the conglomerate of all closure operators on $\mathcal{C}$ with respect to $\mathcal{M}$ ordered as follows: 
$c\leq c'$ if $c_{X}(m)\leq c'_{X}(m)$ for all $m\in \;$sub$X$ and $X\in \mathcal{C}$.
\begin{defn}
	A closure operator $c$ on $\mathcal{C}$ is $idempotent$ if   $c_{X}(c_{X}(m)) = c_{X}(m)$ for all $m\in\;$sub$X$ and $X\in \mathcal{C}$.
\end{defn}
ICL$(\mathcal{C}, \mathcal{M})$ will denote the conglomerate of all idempotent closure operators on $\mathcal{C}$.
\begin{defn} \label{D13}
\cite{MR0286814}	$A\;topogenous\;order$ $\sqsubset$ on $\mathcal{C}$ is a family  $\{\sqsubset_{X}\; |\;X\in \mathcal{C}\}$ of relations, each
	$\sqsubset_{X}$ on sub$X$, such that:
	\begin{itemize}
		\item [$(T1)$] $m\sqsubset_{X} n\Rightarrow m\leq n$ for every $m, n\in \;$sub$X$,
		\item [$(T2)$] $m\leq n\sqsubset_{X} p\leq q\Rightarrow m\sqsubset_{X} q$ for every $m, n, p, q\in \;$sub$X$, and
		\item [$(T3)$] every morphism $f: X\longrightarrow Y$ in $\mathcal{C}$ is $\sqsubset$-continuous, $m\sqsubset_{Y} n\Rightarrow f^{-1}(m)\sqsubset_{X} f^{-1}(n)$
		for every $m, n\in \;$sub$Y$.
	\end{itemize}
\end{defn}
	Given two topogenous orders $\sqsubset$ and $\sqsubset'$ on $\mathcal{C}$, $\sqsubset\subseteq \sqsubset'$ if and only if
$m\sqsubset_{X} n\Rightarrow m\sqsubset'_{X} n$ for all $ m, n\in \;$sub$X$. The resulting ordered congolomerate of all topogenous orders on $\mathcal{C}$ is denoted by TORD$(\mathcal{C}, \mathcal{M}).$\\ 
A topogenous order $\sqsubset$ is said to be 
\begin{itemize}
	\item [$(1)$]$\bigwedge$-$preserving$ if $(\forall i\in I: m\sqsubset_{X} n_{i})\Rightarrow  m\sqsubset_{X} \bigwedge  n_{i}$, and 
	\item [$(2)$] $interpolative$ if $m\sqsubset_{X} n\Rightarrow (\exists\;p)\;|\;m\sqsubset_{X} p\sqsubset_{X} n$ for all $X\in \mathcal{C}$.
\end{itemize}
The ordered conglomerate of all $\bigwedge$-preserving and interpolative topogenous orders is denoted by $\bigwedge$-TORD$(\mathcal{C}, \mathcal{M})$ and INTORD$(\mathcal{C}, \mathcal{M}).$
respectively. $\bigwedge$-INTORD$(\mathcal{C}, \mathcal{M})$  will denote the conglomerate of all interpolative meet preserving topogenous orders. 
\begin{prop}\label{P20}
\cite{MR0286814} $\bigwedge-TORD(\mathcal{C}, \mathcal{M})$ is order isomorphic to $CL(\mathcal{C}, \mathcal{M})$. The inverse assignments of
	each other are given by
	$$c^{\sqsubset}_{X}(m) = \bigwedge\{p\; |\; m\sqsubset_{X} p\}\;\mbox{and}\;\; m\sqsubset^{c}_{X} n\Leftrightarrow c_{X}(m)\leq n\;\mbox{for all}\;X\in \mathcal{C}.$$
\end{prop}
\begin{cor}\label{C21}
	$\bigwedge-INTORD(\mathcal{C}, \mathcal{M})$ is order isomorphic to $ICL(\mathcal{C}, \mathcal{M}).$
\end{cor}
\section{The quasi-uniform structures}
It is well known (see e.g \cite{dowker1962mappings}) that an (entourage) quasi-uniformity on a set $X$ can be equivalently expressed as an appropriate family of maps $U : X\longrightarrow \mathcal{P}(X)$. Since these maps can easily be extended to endomaps on $\mathcal{P}(X)$, it is possible to think of a quasi-uniformity on $\mathcal{C}$ as a suitable family of endomaps on sub$X$ for each $X\in \mathcal{C}$. This is the point expressed in Definition \ref{D21}. Let us denote by  $\mathcal{F}$(sub$X$) the endofunctor category on sub$X$ for each $X\in \mathcal{C}.$ It is clear that for all $U,V \in\mathcal{F}$(sub$X)$, $U\leq V$ if $U(m)\leq V(m)$ for all $m\in $ sub$X$.
\begin{defn}\label{D21}
	\cite{iragi2019quasi} $A\;quasi$-$uniformity$ on $\mathcal{C}$ with respect to $\mathcal{M}$ is a family $\mathcal{U} = \{\mathcal{U}_{X}\;|\;X\in \mathcal{C}\}$ with $\mathcal{U}_{X}$ a full subcategory of $\mathcal{F}($sub$X)$ for each 
	$X$ such that:
	\begin{itemize}
		\item [$(U1)$] For any $U\in \mathcal{U}_{X},\;1_{X}\leq U$,
		\item [$(U2)$] For any $U\in \mathcal{U}_{X}$, there is $U'\in \mathcal{U}_{X}$ such that $U'\circ U'\leq U$,
		\item [$(U3)$] For any $U\in \mathcal{U}_{X}$ and $U\leq U',\;U'\in \mathcal{U}_{X}$,
		\item [$(U4)$] For any $U, U'\in \mathcal{U}_{X}, U\wedge U'\in \mathcal{U}_{X}$,
		\item [$(U5)$] For any $\mathcal{C}$-morphism $f: X\longrightarrow Y$ and $U\in \mathcal{U}_{Y}$, there is $U'\in \mathcal{U}_{X}$ such that $f(U'(m))\leq U(f(m))$ for any $m\in\;$ sub$X$.
	\end{itemize}
\end{defn}
	We shall denote by QUnif$(\mathcal{C}, \mathcal{M})$ the conglomerate of all quasi-uniform structures on $\mathcal{C}$. It is ordered as follows: $\mathcal{U}\leq \mathcal{V}$ if for all $X\in \mathcal{C}$ and $U\in \mathcal{U}_{X}$, there is $V\in \mathcal{V}_{X}$ such that $V\leq U$.
In most cases we describe a quasi-uniformity by defining a base for it. A base for a quasi-uniformity  $\mathcal{U}$ on $\mathcal{C}$ is a family  $\mathcal{B} = \{ \mathcal{B}_{X} \;|X\in \mathcal{C}\}$ with each $\mathcal{B}_{X}$ a full subcategory of $\mathcal{F}$(sub$X$) for all $X\in \mathcal{C}$ satisfying all the axioms in Definition \ref{D21} except $(U3)$. If $\mathcal{B}_{X}$ for any $X\in \mathcal{C}$ is a base element with a single member $V$, we shall write $V_{X}$. A base for quasi-uniformity on $\mathcal{C}$ is $transitive$ if for all $X\in \mathcal{C}$ and $U\in \mathcal{B}_{X}$, $U\circ U = U$. A quasi-uniformity with a transitive base is called $a\;transitive\;quasi$-$uniformity$.
The ordered conglomerate of all transitive quasi-uniformities on $\mathcal{C}$ will be denoted by TQUnif$(\mathcal{C}, \mathcal{M})$. 

\begin{defn}\label{D1}
	\cite{iragi2019quasi}	$A\;syntopogenous\;structure$ on $\mathcal{C}$ with respect to $\mathcal{M}$ is a family \\ $\mathcal{S} = \{\mathcal{S}_{X}\;|\;X\in \mathcal{C}\}$ such that each $\mathcal{S}_{X}$ is a set of relations on sub$X$
	satisfying:
	\begin{itemize}
		\item [$(S1)$] Each $\sqsubset_{X}\in \mathcal{S}_{X}$ is a relation on sub$X$ satisfying $(T1)$ and $(T2),$
		\item [$(S2)$] $\mathcal{S}_{X}$ is a directed set with respect to inclusion,
		\item [$(S3)$] $\sqsubset_{X} = \bigcup \mathcal{S}_{X}$ is an interpolative topogenous order.
	\end{itemize}
\end{defn}
The ordering of topogenous orders can be extended to syntopogenous structures in the following way: $\mathcal{S}\leq \mathcal{S'}$ if for all $X\in \mathcal{C}$ and $\sqsubset_{X}\in \mathcal{S}_{X}$, there is $\sqsubset'_{X}\in \mathcal{S'}_{X}$ such that 
$\sqsubset_{X}\subseteq \sqsubset'_{X}$. The resulting conglomerate will be denoted by SYnt$(\mathcal{C}, \mathcal{M})$. $\mathcal{S}\in $ SYnt$(\mathcal{C}, \mathcal{M})$ is $co$-$perfect$ if each $\sqsubset_{X}\in \mathcal{S}_{X}$ is $\bigwedge$-preserving for all $X\in \mathcal{C}$.
It is $interpolative$ if every $\sqsubset_{X}\in \mathcal{S}_{X}$ interpolates.  The ordered conglomerate of all interpolative co-perfect syntopogenous structures will be denoted by INTCSYnt$(\mathcal{C}, \mathcal{M})$. The ordered conglomerate of all co-perfect syntopogenous structures will be denoted by CSYnt$(\mathcal{C}, \mathcal{M})$. $\mathcal{S}\in $ SYnt$(\mathcal{C}, \mathcal{M})$ is $simple$ if $\mathcal{S}_{X} = \{\sqsubset_{X}\}$ where $\sqsubset_{X}$ is an interpolative topogenous order for any $X\in \mathcal{C}$.

\thm\label{T1}
	\cite{iragi2019quasi}	QUnif$(\mathcal{C}, \mathcal{M})$ is order isomorphic to CSYnt$(\mathcal{C}, \mathcal{M})$.
	The inverse assignments of
	each other $\mathcal{U}\longmapsto \mathcal{S}^{\mathcal{U}}$ and $\mathcal{S}\longmapsto \mathcal{U}^{\mathcal{S}}$  are given by 
	$$\mathcal{S}^{\mathcal{B}}_{X} = \{\sqsubset^{U}_{X}\;|\;U\in \mathcal{B}_{X}\}\;\mbox{where}\;m\sqsubset^{U}_{X} n\Leftrightarrow U(m)\leq n,\;\mbox{and}$$
	$$\mathcal{B}^{\mathcal{S}}_{X} = \{U^{\sqsubset}\;|\;\sqsubset_{X}\in \mathcal{S}_{X}\}\;\mbox{where}\;U^{\sqsubset}(m) = \bigwedge\{n\;|\;m\sqsubset_{X} n\}$$
	for all $X\in \mathcal{C}$ and $m,\;n\in \;$sub$X$.
\endthm

Since $\mathcal{S}_{X}\subseteq \bigwedge$-TORD $(\mathcal{C}, \mathcal{M})$ for each $\mathcal{S}\in$  CSYnt$(\mathcal{C}, \mathcal{M})$, it follows from the above theorem and Proposition \ref{P20} that a quasi-uniformity on $\mathcal{C}$ is a collection of families of closure operators.

By Corollary \ref{C21} ( see also \cite{MR3641225}, Corollary $4.2.3)$, $\bigwedge$-INTORD$(\mathcal{C}, \mathcal{M})$ is isomorphic to the conglomerate of idempotent closure operators and from Theorem \ref{T1}, 
CSYnt$(\mathcal{C}, \mathcal{M})\cong $ QUnif$(\mathcal{C}, \mathcal{M})$. Thus every idempotent closure operator on $\mathcal{C}$ is a quasi-uniformity.
\section{Quasi-uniform structures induced by (co)pointed endofunctors} 

{ \bf{Throughout this section, the class $\mathcal{E}$ will be assumed to be stable under pullbacks along $\mathcal{M}$-morphisms.} }

{\indent	Already the axiom $(S3)$ of Definition \ref{D1} includes the fact that every morphism in $\mathcal{C}$ must be continuous with respect to the syntopogenous structure. In the next definition, we introduce the continuity of $\mathcal{C}$-morphisms with respect to two syntopogenous structures on $\mathcal{C}$. Our aim being to use this definition to construct new syntopogenous structures from old. In particular new quasi-uniformities and new closure operators from old. These are particularly important as they turn
	out to describe initial structures induced by certain maps in spaces.
	
	\begin{defn}\label{D2}
		Let $\mathcal{S},\;\mathcal{S'}\in $ SYnt$(\mathcal{C}, \mathcal{M})$. A morphism $f: X\longrightarrow Y$ is $(\mathcal{S}, \mathcal{S'})$-continous if for all $\sqsubset'_{Y}\in \mathcal{S'}_{Y}$, there is $\sqsubset_{X}\in \mathcal{S}_{X}$
		such that $f(m)\sqsubset'_{Y} n\Rightarrow m\sqsubset_{X} f^{-1}(n)$ for all $m\in \;$sub$X$ and $n\in \;$sub$Y$, equivalently $m\sqsubset'_{Y} n\Rightarrow f^{-1}(m)\sqsubset_{X} f^{-1}(n)$  for all $n, m\in \;$sub$Y.$ 
	\end{defn}
	Since every $\mathcal{C}$-morphism $f$ is $(\mathcal{S}, \mathcal{S})$-continuous and $(\mathcal{S'},\mathcal{S'})$-continuous,  
	$f$ is $(\mathcal{S}, \mathcal{S'})$-continuous if $\mathcal{S'}\leq \mathcal{S}$.
	Because $\mathcal{S}$ is simple if each $\mathcal{S}_{X} = \{\sqsubset_{X}\}$ where $\sqsubset_{X}$ is an interpolative topogenous order, we obtain the following proposition.
	\begin{prop}\label{P512}
		Let $\mathcal{S}\;\mbox{and}\; \mathcal{S'}$ be simple syntopogenous structures i.e $\mathcal{S}_{X} = \{\sqsubset_{X}\}, \mathcal{S'}_{X} = \{\sqsubset'_{X}\}\subseteq $INTORD$(\mathcal{C}, \mathcal{M})$. Then $f$ is $(\mathcal{S}, \mathcal{S'})$-continuous if and only if
		$f(m)\sqsubset'_{Y} n\Rightarrow m\sqsubset_{X} f^{-1}(n)$ for all $m\in \;$sub$X$ and $n\in \;$sub$Y$.
	\end{prop}
	The next proposition is obtained from Theorem \ref{T1}.
	\begin{prop}
		If $\mathcal{S},\;\mathcal{S'}\in $ SYnt$(\mathcal{C}, \mathcal{M})$. Then $f$ is $(\mathcal{S}, \mathcal{S'})$-continuous if and only if for any $V\in \mathcal{B}^{\mathcal{S'}}_{Y}$ there is $U\in \mathcal{B}^{\mathcal{S}}_{X}$ 
		such that $f(U(m))\leq V(f(m))$ for all $m\in \;$sub$X$.
	\end{prop}
	\begin{proof}
		Assume that $f: X\longrightarrow Y$ is $(\mathcal{S}, \mathcal{S'})$-continuous and $\mathcal{S},\;\mathcal{S'}\in $ SYnt$(\mathcal{C}, \mathcal{M})$. Then for any 
		$V\in \mathcal{B}^{\mathcal{S'}}_{Y}$, there is $\sqsubset'_{Y}\in \mathcal{S'}_{Y}$ which determines $V$ and there is $\sqsubset_{X}\in \mathcal{S}_{X}$ 
		such that $f(m)\sqsubset'_{Y} n\Rightarrow m\sqsubset_{X} f^{-1}(n)$.
		Now $U(m) = U^{\sqsubset}_{X}(m) = \bigwedge\{p\;|\;m\sqsubset_{X} p\}\leq \bigwedge \{f^{-1}(n)\;|\;f(m)\sqsubset'_{Y} n\} = f^{-1}(V(f(m)))\Rightarrow
		U(m)\leq f^{-1}(V(f(m))\Leftrightarrow f(U(m))\leq V(f(m)).$ Conversely, assume that for any $V\in \mathcal{B}^{\mathcal{S}}_{Y}$ there is $U\in \mathcal{B}^{\mathcal{S}}_{X}$ 
		such that $f(U(m))\leq V(f(m))$.  Now, for any $\sqsubset'_{Y}\in \mathcal{S'}_{Y}$, there is, by Theorem \ref{T1}, $V\in \mathcal{B}^{\mathcal{S}}$ such that $\sqsubset_{Y} = \sqsubset^{V}$.
		Thus $f(m)\sqsubset'_{Y} n\Leftrightarrow V(f(m))\leq n\Rightarrow f(U(m))\leq n\Leftrightarrow U(m)\leq f^{-1}(n)\Leftrightarrow m\sqsubset^{U}_{X} f^{-1}(n)\Leftrightarrow m\sqsubset_{X} f^{-1}(n).$
	\end{proof}
	The proposition above provides us with the next definition.
	\begin{defn}
		Let $\mathcal{U},\; \mathcal{U'}\in $ QUnif$(\mathcal{C}, \mathcal{M})$ and $f : X\longrightarrow Y$ a $\mathcal{C}$-morphism. $f$ is $(\mathcal{U},\; \mathcal{U'})$-continous if for any $U'\in \mathcal{U'}_{Y}$, there is $U\in \mathcal{U}_{X}$ such that $f(U(m))\leq U'(f(m))$ for all $m\in \;$sub$X$. 
	\end{defn}
	
	Propositions \ref{P512} and Corollary \ref{C21} allow us to prove the following.
	\begin{prop}
		Let $\mathcal{S}\;\mbox{and}\; \mathcal{S'}$ be simple and co-perfect syntopogenous structures i.e $\mathcal{S}_{X} = \{\sqsubset_{X}\}, \mathcal{S'}_{X} = \{\sqsubset'_{X}\}\subseteq \bigwedge-INTORD(\mathcal{C}, \mathcal{M})$. Then $f$ is $(\mathcal{S}, \mathcal{S'})$-continuous if and only if
		$f(c^{\sqsubset}_{X}(m))\leq c^{\sqsubset'}_{X}(f(m))$ for all $m\in \;$sub$X$.
	\end{prop}
	\begin{defn}
		\cite{MR2122803} Let $c,c'\in $ CL$(\mathcal{C}, \mathcal{M})$ and $f : X\longrightarrow Y$ a $\mathcal{C}$-morphism. $f$ is $(c, \;c')$-continuous if $f(c_{X}(m))\leq c'_{X}(f(m))$ for all $m\in \;$sub$X$.
	\end{defn}
	For a syntopogenous structure $\mathcal{S}$ on $\mathcal{C}$ and a class $\mathcal{F}$ of $\mathcal{C}$-morphisms, we ask if there is a coarsest syntopogenous structure $\mathcal{S'}$ on $\mathcal{C}$ for which every morphism in $\mathcal{F}$ is $(\mathcal{S'}, \mathcal{S})$-continous. In the next theorem, we provide an answer to this question in the case $\mathcal{F} = \{\eta_{X}:X\in \mathcal{C}\}$, for an $\mathcal{E}$-pointed endofunctor $(F, \eta)$ of $\mathcal{C}$. Later on we shall deal with a somehow  
	dual case. Let us also note that a similar question has been  asked in the case of a closure operator (see \cite{MR2122803}, chapter 5). We prove that the results obtained in (\cite{MR2122803}) can be deduced from those we prove here. 
	\thm\label{T222}
		Let $(F,\eta)$ be an $\mathcal{E}$-pointed endofunctor of $\mathcal{C}$ and $\mathcal{S}$ a syntopogenous structure on $\mathcal{C}$ with respect to $\mathcal{M}$.  Then
		$\mathcal{S}^{F, \eta}_{X} = \{\sqsubset^{F, \eta}_{X}\;|\;\sqsubset_{FX}\in \mathcal{S}_{FX}\}$
		with $m\sqsubset_{X}^{F, \eta} n\Leftrightarrow \eta_{X}(m)\sqsubset_{FX} p$ and $\eta_{X}^{-1}(p)\leq n$ for some $p\in \mbox{sub}FX$ is the coarsest syntopogenous structure on $\mathcal{C}$ with respect to $\mathcal{M}$ for which 
		every $\eta_{X} : X\longrightarrow FX$ is $(\mathcal{S}^{F,\eta}, \mathcal{S})$-continuous.
		If $\mathcal{S}$ is interpolative (co-perfect), then $\mathcal{S}^{F,\eta}$ is interpolative (co-perfect, respectively).
	\endthm
	\begin{proof}
		$\mathcal{S}^{F, \eta}$ is clearly a syntopogenous structure and $\eta_{X}$ is $(\mathcal{S}^{F,\eta}, \mathcal{S})$-continuous, since for all $\sqsubset_{X}\in \mathcal{S}_{X}, 
		\eta_{X}(m)\sqsubset_{FX} n\Rightarrow \eta_{X}(m)\sqsubset_{FX} (\eta_{X}(\eta_{X}^{-1}(n))\Leftrightarrow m\sqsubset_{X}^{F,\eta} \eta_{X}^{-1}(n).$\\                                                                           
		If $\mathcal{S'}$ is another syntopogenous structure on $\mathcal{C}$ such that $\eta_{X}$ is $(\mathcal{S'},\;\mathcal{S})$-continuous, then for any $\sqsubset^{F,\eta}_{X}\in \mathcal{S}^{F,\eta}_{X}$, 
		$m\sqsubset^{F,\eta}_{X} n\Leftrightarrow \eta_{X}(m)\sqsubset_{FX} p$ and $\eta^{-1}_{X}(p)\leq n$. This implies that there is $\sqsubset'_{X}\in \mathcal{S'}_{X}$
		such that $m\sqsubset'_{X} \eta_{X}^{-1}(p)\leq n \Rightarrow m\sqsubset'_{X} n.$ Thus $\mathcal{S}^{F,\eta}\leq \mathcal{S'}$.
		
		If $\mathcal{S}$ is interpolative and $m\sqsubset_{X}^{F, \eta} n$, then $\eta_{X}(m)\sqsubset_{FX} p$ and $\eta_{X}^{-1}(p)\leq n$ for some $p\in \;$sub$FX$. This implies that there is $l\in \;$sub$FX$
		such that $\eta_{X}(m)\sqsubset_{FX} l \sqsubset_{FX} p$. Thus $\eta_{X}(m)\sqsubset_{FX} \eta_{X}(\eta^{-1}_{X}(l)) \sqsubset_{FX} p$, that is 
		$m\sqsubset_{X}^{F, \eta} \eta_{X}^{-1}(l)\sqsubset_{X}^{F, \eta}n.$ It is also not hard to see that $\mathcal{S}^{F,\eta}$ is co-perfect if $\mathcal{S}$ has the same property.
		
	\end{proof}
	Viewing a reflector as endofunctor of $\mathcal{C}$, one obtains the proposition below.
	\begin{cor}\label{C222}
		Let $\mathcal{A}$ be $\mathcal{E}$-reflective subcategory of $\mathcal{C}$ and $\mathcal{S}$ a syntopogenous structure on $\mathcal{A}$ with respect to $\mathcal{M}$.  Then
		$\mathcal{S}^{\mathcal{A}}_{X} = \{\sqsubset^{\mathcal{A}}_{FX}\;|\;\sqsubset_{FX}\in \mathcal{S}_{FX}\}$
		with $m\sqsubset_{X}^{\mathcal{A}} n\Leftrightarrow \eta_{X}(m)\sqsubset_{FX} p$ and $\eta_{X}^{-1}(p)\leq n$ for some $p\in \mbox{sub}FX$  is the coarsest syntopogenous structure on $\mathcal{C}$ with respect to $\mathcal{M}$ for which 
		every reflection morphism $\eta_{X} : X\longrightarrow FX$ is $(\mathcal{S}^{\mathcal{A}}, \mathcal{S})$-continous.
		If $\mathcal{S}$ is interpolative (co-perfect), then $\mathcal{S}^{\mathcal{A}}$ is interpolative (co-perfect, respectively).
	\end{cor}
	Since $\mathcal{S}^{F,\eta}$ is co-perfect provided $\mathcal{S}$ is co-perfect, Theorem \ref{T1} gives us the next proposition. 
	\begin{prop}
		Let $(F,\eta)$ be a pointed endofunctor of $\mathcal{C}$ and $\mathcal{S}\in $CSYnt$(\mathcal{C}, \mathcal{M})$.  Then
		$$\mathcal{B}^{\mathcal{S}^{F,\eta}}_{X} = \{U^{\sqsubset^{F,\eta}}\;|\;U^{\sqsubset}\in \mathcal{B}^{\mathcal{S}}_{FX}\}\;\mbox{with}\; U^{\sqsubset^{F,\eta}}(m) = \eta_{X}^{-1}(U^{\sqsubset}(\eta_{X}(m)))$$ is a base for the coarsest quasi-uniformity on
		$\mathcal{C}$ with respect to $\mathcal{M}$ for which every \\ $\eta_{X}:X\longrightarrow FX$ is 
		$(\mathcal{U}^{\mathcal{S}^{F,\eta}}, \mathcal{U}^{\mathcal{S}})$-continous. $\mathcal{B}^{\mathcal{S}^{F,\eta}}$ is a transitive base provided that $\mathcal{S}$ is interpolative.
	\end{prop}
	
	\begin{proof}
		$(U1)$, $(U2)$ and $(U4)$ are clear.
		$(U5)$ Let $f : X\longrightarrow Y$ be a $\mathcal{C}$-morphism and  $U^{\sqsubset^{F,\eta}}\in \mathcal{B}^{\mathcal{S}^{F,\eta}}_{Y}$ for $\sqsubset_{FY}\in \mathcal{S}_{FY}$. Then there is $\sqsubset_{FX}\in \mathcal{S}_{FX}$ such that 
		$f(V^{\sqsubset_{FX}}(m))\leq U^{\sqsubset_{FY}}(f(m))$.
		\begin{align*}
		\mbox{Thus}\;\;f(V^{\sqsubset^{F,\eta}}(m)) & = f(\eta^{-1}_{X}(V^{\sqsubset_{FX}}(\eta_{X}(m)))\\ 
		& \leq \eta_{Y}^{-1}(Ff)(V^{\sqsubset_{FX}}(\eta_{X}(m))) &&\text{Lemma \ref{L22}}\\
		& \leq \eta_{Y}^{-1}(U^{\sqsubset_{FY}}(Ff)(\eta_{X}(m)))\\
		& = \eta^{-1}_{Y}(U^{\sqsubset_{FX}}(\eta_{Y}(f(m))) &&\text{ Definition \ref{D154}}\\
		& = U^{\sqsubset^{F, \eta}}(f(m))\\                                             
		\end{align*}
		Since, for any $\sqsubset_{X}\in \mathcal{S}_{X}$, $U^{\sqsubset^{F,\eta}}(m) = \eta_{X}^{-1}(U^{\sqsubset}(\eta_{X}(m)))\Rightarrow \eta_{X}(U^{\sqsubset^{F,\eta}}(m))\leq U^{\sqsubset}(\eta_{X}(m))$, 
		$\eta_{X}$ is $(\mathcal{U}^{\mathcal{S}^{F,\eta}},\;\mathcal{U})$-continous for all $X\in \mathcal{C}$.
		If $\mathcal{S}$ is interpolative  then 
		$U^{\sqsubset^{F,\eta}}(U^{\sqsubset^{F,\eta}}(m)) = U^{\sqsubset^{F,\eta}}(\eta^{-1}_{X}(U^{\sqsubset}(\eta_{X}(m)))
		= \eta^{-1}_{X}(U^{\sqsubset}(\eta_{X}(\eta^{-1}_{X}(U^{\sqsubset}(\eta_{X}(m))))
		\leq \eta^{-1}_{X}(U^{\sqsubset}(U^{\sqsubset}(\eta_{X}(m)))) 
		=\\ \eta^{-1}_{X}(U^{\sqsubset}(\eta_{X}(m))) 
		= U^{\sqsubset^{F,\eta}}(m).$\\
		Let $\mathcal{B'}$ be a base for another quasi-uniformity $\mathcal{U'}$ on $\mathcal{C}$ such that $\eta_{X}$ is $(\mathcal{U'},\;\mathcal{U}^{\mathcal{S}})$-continuous, then 
		for any $U^{\sqsubset}\in \mathcal{B}^{\mathcal{S}}_{FX}$, there is $U'\in \mathcal{B'}_{X}$ such that $\eta_{X}(U'(m))\leq U^{\sqsubset}(\eta_{X}(m))\Leftrightarrow
		U'(m)\leq \eta_{X}^{-1}(U^{\sqsubset}(\eta_{X}(m))) = U^{\sqsubset^{F,\eta}}(m)$. Thus $\mathcal{B}^{\mathcal{S}^{F, \eta}}\leq \mathcal{B'}$.  
	\end{proof}
	One sees from the proof of the above proposition that the condition of $(F,\eta)$ being $\mathcal{E}$-pointed is not needed when the syntopogenous structure is co-perfect.
	\begin{prop}\label{C8}
		Let $(F,\eta)$ be a pointed endofunctor of $\mathcal{C}$ and $\mathcal{S}$ be simple and co-perfect syntopogenous structures i.e $\mathcal{S}_{X} = \{\sqsubseteq_{X}\}\in \bigwedge- $INTORD$(\mathcal{C}, \mathcal{M})$. 
		Then $c^{\sqsubset^{F,\eta}}(m) = \\ \eta_{X}^{-1}(c^{\sqsubset}_{FX}(\eta_{X}(m)))$  is an idempotent closure operator. It is the largest closure operator 
		on $\mathcal{C}$ for which every $\eta_{X} : X\longrightarrow FX$ is  $(c^{\sqsubset^{F,\eta}}, c^{\sqsubset})$-continuous.
	\end{prop}
	The above closure operator was first introduced on the category of topological spaces and continuous maps by L. Stramaccia (\cite{MR958752}), then on topological categories by D. Dikranjan (\cite{dikranjan1992semiregular}) and later on an arbitrary category by Dikranjan and Tholen (\cite{MR2122803}). It is a special case of the pullback closure studied by D. Holgate in \cite{holgate1996pullback, holgate1995pullback}.
	\begin{cor}\label{C1}
		Let $\mathcal{A}$ be a reflective subcategory of $\mathcal{C}$  and $\mathcal{S}$ a co-perfect syntopogenous structure on $\mathcal{A}$ with respect to $\mathcal{M}$.  Then
		$$\mathcal{B}^{\mathcal{A}}_{X} = \{U^{\sqsubset^{\mathcal{A}}}\;|\;U^{\sqsubset}\in \mathcal{B}^{\mathcal{S}}_{FX}\}\;\mbox{with}\; U^{\sqsubset^{\mathcal{A}}}(m) = \eta_{X}^{-1}(U^{\sqsubset}(\eta_{X}(m)))$$ is a base for the coarsest quasi-uniformity on
		$\mathcal{C}$ with respect to $\mathcal{M}$ for which every reflection morphism $\eta_{X}:X\longrightarrow FX$ is 
		$(\mathcal{U}^{\mathcal{S}^{\mathcal{A}}},\;\mathcal{U}^{\mathcal{S}})$-continous. $\mathcal{B}^{\mathcal{S}^{\mathcal{A}}}$ is a transitive base provided that $\mathcal{S}^{F,\eta}$ is interpolative.
	\end{cor}
	
	Corollary \ref{C1}  allows us to obtain the quasi-uniform structure induced by any reflective subcategory of $\bf QUnif$ and to conclude that it is the initial quasi-uniformity for which the reflection map is quasi-uniformly continous (see Example $6.1$). 
	
	\thm
		Let $(G,\varepsilon)$ be a $\mathcal{M}$-copointed endofunctor of $\mathcal{C}$ and $\mathcal{S}$ a syntopogenous structure on $\mathcal{C}$, then  
		$\mathcal{S}^{G, \varepsilon}_{X} = \{\sqsubset^{G,\varepsilon}_{X}\;|\;\sqsubset_{GX}\in \mathcal{S}_{GX}\}$ with $m\sqsubset^{G, \varepsilon}_{X} n\Leftrightarrow \varepsilon^{-1}_{X}(n)\sqsubset_{GX} \varepsilon^{-1}_{X}(n)$ for all $m\in \;$sub$X$ and $n\geq m$, is the finest syntopogenous structure on 
		$\mathcal{C}$ for which every $\varepsilon_{X} : GX\longrightarrow X$ is $(\mathcal{S}, \mathcal{S}^{G, \varepsilon})$-continuous. 
	\endthm
	\begin{proof}
		A routine check shows that $\mathcal{S}^{G, \varepsilon}$ is a syntopogenous structure on $\mathcal{C}$.
		For all $X\in \mathcal{C}$, $\varepsilon_{X} : GX\longrightarrow X$ is $(\mathcal{S},\;\mathcal{S}^{G,\varepsilon})$-continuous, since for any 
		$\sqsubset^{G,\varepsilon}_{X}\in \mathcal{S}^{G,\varepsilon}_{X}$ and $m,n\in \;$sub$X$ with $n\leq m$, $m\sqsubset^{G,\varepsilon}_{X}n\Rightarrow \varepsilon^{-1}_{X}(n)\sqsubset_{GX} \varepsilon^{-1}_{X}(n)$.
		
		If $\mathcal{S'}$ is another syntopgenous structure on $\mathcal{C}$ such that 
		$\varepsilon_{X}$ is $(\mathcal{S}, \mathcal{S'})$-continuous, then for any $\sqsubset_{X}\in\mathcal{S'}_{X}$, $m\sqsubset'_{X} n\Rightarrow \varepsilon_{X}(\varepsilon_{X}^{-1}(m))\sqsubset'_{X} n \Rightarrow
		\exists\sqsubset_{GX}\in \mathcal{S}_{GX}\;|\; \varepsilon_{X}^{-1}(m)\sqsubset_{X} \varepsilon_{X}^{-1}(n)\Leftrightarrow m\sqsubset^{G,\varepsilon}_{X} n$.
	\end{proof}
	\begin{cor}
		Let $\mathcal{A}$ be an $\mathcal{M}$-coreflective subcategory of $\mathcal{C}$ and $\mathcal{S}$ a syntopogenous structure on $\mathcal{A}$, then   
		$\mathcal{S}^{\mathcal{A}}_{X} = \{\sqsubset^{\mathcal{A}}_{X}\;|\;\sqsubset_{GX}\in \mathcal{S}_{GX}\}$ with $m\sqsubset^{\mathcal{A}}_{X} n\Leftrightarrow \varepsilon^{-1}_{X}(n)\sqsubset_{GX} \varepsilon^{-1}_{X}(n)$ for all $m\in \;$sub$X$ and $n\geq m$, is the finest syntopogenous structure on 
		$\mathcal{C}$ for which every coreflection $\varepsilon_{X} : GX\longrightarrow X$ is $(\mathcal{S}, \mathcal{S}^{\mathcal{A}})$-continuous.
	\end{cor}

	\begin{prop}
		Assume that $f^{-1}$ commutes with the join of subobjects for any $f\in \mathcal{C}$. Let $(G,\varepsilon)$ be an $\mathcal{M}$-copointed endofunctor of $\mathcal{C}$ and $\mathcal{S}\in $CSYnt$(\mathcal{C}, \mathcal{M})$.  Then
		$$\mathcal{B}^{\mathcal{S}^{G,\varepsilon}}_{X} = \{V^{\sqsubset^{G,\varepsilon}}\;|\;V^{\sqsubset}\in \mathcal{B}^{\mathcal{S}}_{GX}\}\;\mbox{with}\; V^{\sqsubset^{F,\varepsilon}}_{X}(m) = m\vee \varepsilon_{X}(V^{\sqsubset}(\varepsilon_{X}^{-1}(m)))$$
		is a base for the finest quasi-uniformity on $\mathcal{C}$ which makes every $\varepsilon_{X}$ $(\mathcal{V}, \mathcal{V}^{G,\varepsilon})$-continous.\\
		
	\end{prop}
	\begin{proof}
		It is not hard to check that $\mathcal{B}^{\mathcal{S}^{G,\varepsilon}}_{X}$ is a base for a quasi-uniformity on $\mathcal{C}$.
		Since $\varepsilon_{X}(V^{\sqsubset}(\varepsilon_{X}^{-1}(m)))\leq V^{\sqsubset^{G,\varepsilon}}(m)\Leftrightarrow V^{\sqsubset}(\varepsilon_{X}^{-1}(m))\leq \varepsilon_{X}^{-1}(V^{\sqsubset^{G,\varepsilon}}(m)), \varepsilon_{X}$ is 
		$(\mathcal{V}, \mathcal{V}^{\mathcal{S}^{G,\varepsilon}})$-continous.\\
		Let $\mathcal{B'}$ be base for another quasi-uniformity $\mathcal{V'}$ on $\mathcal{C}$ such that $\varepsilon_{X}$ is $(\mathcal{V}, \mathcal{V'})$-continuous. Then for all $V'\in \mathcal{V'}_{X}$, there is $V\in \mathcal{V}_{GX}$
		such $V(\varepsilon^{-1}_{X}(m))\leq \varepsilon^{-1}_{X}(V'(m))\Leftrightarrow  \varepsilon_{X}(V(\varepsilon^{-1}_{X}(m)))\\
		\leq V'(m)\Rightarrow m\vee \varepsilon_{X}(V(\varepsilon^{-1}_{X}(m)))\leq V'(m)\Leftrightarrow V^{\sqsubset^{G,\varepsilon}}(m)\leq V'(m)$. 
		Thus $\mathcal{B'}\leq \mathcal{B}^{G,\varepsilon}.$
	\end{proof}
	\begin{prop}
		Let $(G,\varepsilon)$ be a copointed endofunctor of $\mathcal{C}$ and $\mathcal{S}$ be simple and co-perfect syntopogenous structure i.e $\mathcal{S}_{X} = \{\sqsubset_{X}\}\in \bigwedge-$INTORD$(\mathcal{C}, \mathcal{M})$, then  for all $m\in \;$sub$X$, 
		$c^{\sqsubset^{G, \varepsilon}}(m) = m\vee \varepsilon_{X}(c^{\sqsubset}_{GX}(\varepsilon^{-1}_{X}(m)))$ is is an idempotent closure operator on $\mathcal{C}$. 
		It is the least closure operator for which every $\varepsilon_{X} : GX\longrightarrow X$ is $(c, c^{G, \varepsilon})$-continuous. 
	\end{prop}
	\begin{cor}\label{C2}
		Assume that $f^{-1}$ commutes with the join of subobjects for any $f\in \mathcal{C}$. Let $\mathcal{A}$ be an $\mathcal{M}$-coreflective subcategory of $\mathcal{C}$ and $\mathcal{S}$ a syntopogenous $\mathcal{A}$.  Then
		$$\mathcal{B}^{\mathcal{A}}_{X} = \{V^{\sqsubset^{\mathcal{A}}}\;|\;V^{\sqsubset}\in \mathcal{B}^{\mathcal{S}}_{GX}\}\;\mbox{with}\; V^{\sqsubset^{\mathcal{A}}}(m) = m\vee \varepsilon_{X}(V^{\sqsubset}(\varepsilon_{X}^{-1}(m))$$
		is a base for finest quasi-uniformity on $\mathcal{C}$ which makes every coreflection morphism $\varepsilon_{X}$ $(\mathcal{V}, \mathcal{V}^{\mathcal{A}})$-continous. 
	\end{cor}
\section{The continuity of functors with respect to quasi-uniform structures}
Let $\mathcal{A}$ be a category endowed with an $(\mathcal{E'}, \mathcal{M'})$-factorization system for morphisms
and  $\mathcal{A}$ be $\mathcal{M'}$-complete.
	\begin{defn}\cite{MR2122803}
		A functor $F : \mathcal{A}\longrightarrow \mathcal{C}$ is said to preserve subobjects provided that $Fm$
		is an $\mathcal{M}$-subobject for every $\mathcal{M'}$-subobject $m$. It preserves inverse images (resp.  images) of subobjects if $Ff^{-1}(n) = (Ff)^{-1}(Fn)$ (resp. 
		$(Ff)(Fm) = F(f(m))$) for any $\mathcal{A}$-morphism 
		$f : X\longrightarrow Y$ and subobjects $n\in \;$sub$Y$,     $m\in \;$sub$X$.
	\end{defn}
	\begin{defn}\label{D35}
		Let $F : \mathcal{A}\longrightarrow \mathcal{C}$ be a functor that preserves subobjects, $\mathcal{U}\in $\\ QUnif$(\mathcal{A}, \mathcal{M'})$ and $\mathcal{V}\in  $ QUnif$(\mathcal{C}, \mathcal{M})$. Then
		$F$ is $(\mathcal{U},\;\mathcal{V})$-continuous if for all $V\in \mathcal{V}_{FX}$, there is $U\in \mathcal{U}_{X}$ such that 
		$FU(m)\leq V(Fm)$ for all $m\in \;$sub$X$, $X\in \mathcal{A}$.
	\end{defn}
	It can be easily seen that our definition for $(\mathcal{U},\;\mathcal{V})$-continuity of $F$ is a  generalization of $\mathcal{U}$-continuity
	of morphisms to functors. Using Theorem \ref{T1}, we can formulate an equivalent definition of the $(\mathcal{U},\;\mathcal{V})$-continuity of $F$ in terms of co-perfect syntopogenous 
	structures so that $F$ is $(\mathcal{S}, \mathcal{S'})$-continuous will mean that $F$ is continuous with respect to the quasi-uniform structures associated 
	with $\mathcal{S}$ and $\mathcal{S'}$.
	\begin{prop}
		Let $F : \mathcal{A}\longrightarrow \mathcal{C}$ be a functor that preseves subobjects, $\mathcal{S}\in $\\ QUnif$(\mathcal{A}, \mathcal{M'})$ and $\mathcal{S'}\in $ QUnif$(\mathcal{C}, \mathcal{M})$.
		Then $F$ is $(\mathcal{S},\;\mathcal{S'})$-continuous if for all $\sqsubset'_{FX}\in \mathcal{S'}_{FX}$, there is $\sqsubset_{X}\in \mathcal{S}_{X}$ such that 
		$FU^{\sqsubset}(m)\leq U^{\sqsubset'}(Fm)$ for all $m\in \;$sub$X$, $X\in  \mathcal{A}$.
	\end{prop}
	Continuity of a functor between categories supplied with fixed closure operators has been studied in \cite{MR2122803}. We next use the above proposition together with Corollary \ref{C21} and the fact that $\bigwedge-INTORD(\mathcal{C}, \mathcal{M})$ is equivalent to the 
	simple co-perfect syntopogenous structures to produce the $(\mathcal{U},\;\mathcal{V})$-continuity of $F$ in terms of idempotent closure 
	operators.
	\begin{prop}
		Let $F : \mathcal{A}\longrightarrow \mathcal{C}$ be a functor that preseves subobjects, \\ $\mathcal{S}\in $ 
		CSYnt$(\mathcal{A}, \mathcal{M'})$ and $\mathcal{S}\in $ CSYnt$ (\mathcal{C}, \mathcal{M})$ with $\mathcal{S}$ and $\mathcal{S'}$ being simple i.e $\mathcal{S}_{X} = \{\sqsubset_{X}\}$ and 
		$\mathcal{S'}_{FX} = \{\sqsubset'_{FX}\}$. Then $F$ is $(\mathcal{S},\;\mathcal{S'})$-continuous if and only if  for all 
		$Fc^{\sqsubset}_{X}(m)\leq c^{\sqsubset'}_{FX}(Fm)$ for all $m\in \;$sub$X$, $X\in \mathcal{A}$.
	\end{prop}
	
	\begin{defn}
		\cite{MR2122803} Let $F : \mathcal{A}\longrightarrow \mathcal{C}$  a faithful functor. $F$ is called a fibration if every $g : A\longrightarrow FY$ has an $F$-initial ($F$-cartesian) lifting. If we require the existence of an $F$-cartesian lifting of $g : A\longrightarrow FY$ only if $g\in \mathcal{M}$, then $F$ is called an $\mathcal{M}$-fibration.
	\end{defn}
	Let us denote by $IniF$  the class of all $F$-initial morphisms in $\mathcal{A}$.
	Then for an $\mathcal{M}$-fibration $F : \mathcal{A}\longrightarrow \mathcal{C}$, $(\mathcal{E}_{F}, \mathcal{M}_{F})$ where 
	$\mathcal{E}_{F} =  F^{-1}\mathcal{E} = \{e\in \mathcal{C}\;|\;Fe\in \mathcal{E}\}$ and  
	$\mathcal{M}_{F} =  F^{-1}\mathcal{M}\bigcap Ini F$ is a factorization system in $\mathcal{A}$ and
	$\mathcal{M}$-subobject properties in $\mathcal{C}$ are inherited by $\mathcal{M}_{F}$-subobjects in $\mathcal{A}$.

	In particular, 
	\begin{itemize}
		\item [$(1)$] $\mathcal{A}$ has $\mathcal{M}_{F}$-pullbacks if  $\mathcal{C}$ has $\mathcal{M}$-pullbacks.
		\item [$(2)$] $\mathcal{A}$ is $\mathcal{M}_{F}$-complete if $\mathcal{C}$ is $\mathcal{M}$-complete.
		\item [$(3)$] the $\mathcal{M}_{F}$-images and $\mathcal{M}_{F}$-inverse images are obtained by initially lifting $\mathcal{M}$-images and $\mathcal{M}$-inverse images. Consequently 
		$F$ preserves images and inverse images of subobjects. 
	\end{itemize}
	\lem
		\cite{MR2122803} Let $F : \mathcal{A}\longrightarrow \mathcal{C}$ be a faithful $\mathcal{M}$-fibration. 
		\begin{itemize}
			\item [$(1)$]  For any $X\in \mathcal{A}$, sub$X$ and sub$FX$ are order equivalent with the inverse assignments, 
			$\gamma_{X} : subX\longrightarrow subFX$ and $\delta_{X} : subFX\longrightarrow subX$, given 
			by $\gamma_{X}(m) = Fm$ and $\delta_{X}(n) = p$ with $Fp = n$ and $p\in IniF$.
			\item [$(2)$] For any $f : X\longrightarrow Y \in \mathcal{A}$ and suitable subobjects $n, m, n'$ and $m'$.
			\begin{itemize}
				\item [$(1)$] $\gamma_{Y}(f(m)) = (Ff)(\gamma_{X}(m))$.
				\item [$(2)$] $f(\delta_{X}(n)) = \delta_{Y}(Ff)(n)$.
				\item [$(3)$] $f^{-1}(\delta_{Y}(m')) = \delta_{X}((Ff)^{-1}(m'))$.
				\item [$(4)$] $\gamma_{X}(f^{-1}(n')) = (Ff)^{-1}(\gamma_{Y}(n'))$.
			\end{itemize}
		\end{itemize}
	\endlem
	\begin{prop}\label{38}
		Let $F : \mathcal{A}\longrightarrow \mathcal{C}$ be faithful $\mathcal{M}$-fibration and $\mathcal{S}$ be a syntopogenous structure on 
		$\mathcal{C}$ with respect to $\mathcal{M}$. Then 
		$$\mathcal{S}^{F}_{X} = \{\sqsubset^{F}_{X}\;|\;\sqsubset_{FX}\in \mathcal{S}_{FX}\}\;\;\mbox{where}\;\;
		m\sqsubset^{F}_{X} n \Leftrightarrow Fm\sqsubset_{FX} \gamma_{X}(n)$$ is a syntopogenous structure on $\mathcal{A}$ with respect to $\mathcal{M}_{F}$ which is interpolative, co-perfect provided $\mathcal{S}$ has the 
		same properties. Moreover, an $\mathcal{A}$-morphism $f$ is $\mathcal{S}^{F}$-initial provided $Ff$ is $\mathcal{S}$-initial.
	\end{prop}
	\thm\label{257}
		Let $F : \mathcal{A}\longrightarrow \mathcal{C}$ be a faithful $\mathcal{M}$-fibration and $\mathcal{B}$ be a base for a quasi-uniform structure on 
		$\mathcal{C}$ with respect to $\mathcal{M}$. Then $\mathcal{B}^{F}_{X} = \{U^{F}\;|\;U\in \mathcal{B}_{FX}\}\;\;\mbox{where}\;\;U^{F}(m) = \delta_{X}(U(Fm))$ is a base for quasi-uniformity on 
		$\mathcal{A}$ with respect to $\mathcal{M}_{F}$.  It is the coarsest quasi-uniformity for which $F$ is $(\mathcal{U}^{F}, \mathcal{U})$-continuous.
		$\mathcal{B}^{F}$ is transitive provided that $\mathcal{B}$ is a transitive base. Moreover an $\mathcal{A}$-morphism $f$ is $\mathcal{U}^{F}$-initial provided $Ff$ is $\mathcal{S}$-initial.
	\endthm
	\begin{proof}
		It is clear that $\mathcal{B}^{F}$ is a base for a quasi-uniformity on $\mathcal{A}$ which is transitive if $\mathcal{B}$ is transitive. 
		$F$ is $(\mathcal{U}^{F}, \mathcal{U})$-continuous, since for any $U\in \mathcal{B}_{FX}$, $U^{F}(m) = \delta_{X}(U(Fm))
		\Leftrightarrow \gamma_{X}(U^{F}(m)) = U(Fm)\Leftrightarrow F(U^{F}(m)) = U(Fm).$ 
		If $\mathcal{B'}$ is a base for another quasi-uniformity $\mathcal{U'}$ on $\mathcal{A}$ such that $F$ is $(\mathcal{U'}, \mathcal{U})$-continuous, then for all $U^{F}\in \mathcal{B}^{F}_{X}$, there 
		is $U'\in \mathcal{B'}$ such that $FU'(m)\leq U(Fm) = FU^{F}(m)$. Thus  $U'(m) = \delta_{X}(FU'(m))\leq \delta_{X}(FU^{F}(m)) = U^{F}(m)$, that is $\mathcal{B}^{F}\leq \mathcal{B'}$.
		If $Ff$ is $\mathcal{U}$-initial and $U^{F}\in \mathcal{U}^{F}_{X}$, there is $U'\in \mathcal{U}_{FY}$ such that 
		$(Ff)^{-1}(U'(Ff)(p))\leq U(p)$ for all $p\in \;$sub$FX$. Now $ f^{-1}(U'^{F}(f(m))) = f^{-1}(\delta_{Y}(U'(Ff(m))))\\ =  \delta_{X}((Ff)^{-1}(U'(Ff(m)))) =  \delta_{X}((Ff)^{-1}(U'((Ff)(Fm))))\leq \delta_{X}(U(Fm)) = U^{F}(m)$ for all $m\in \;$sub$X$.
	\end{proof}
	
	\begin{cor}
		Under the assumptions of Theorem \ref{257} and $F$ is essentially surjective on objects, then $\mathcal{B}$ is the base of the finest quasi-uniformity
		on $\mathcal{C}$ for which $F$ is $(\mathcal{U}^{F},\mathcal{U})$-continous.
	\end{cor}
	\begin{proof}
		By essential surjectivity of $F$ on objects, we have that for all $Y\in \mathcal{C}$, $Y \cong FX$ for some $X\in \mathcal{A}$. Thus if $\mathcal{B'}$ is another quasi-uniformity on 
		$\mathcal{C}$ such that $F$ is $(\mathcal{U}^{F},\mathcal{U'})$-continuous, then for all $Y\in \mathcal{C}$ and $U'\in \mathcal{U'}_{Y}$, there is $X\in \mathcal{A}$ and $U^{F}\in \mathcal{B}^{F}$ such that 
		$Y \cong FX$ and $FU^{F}(m)\leq U'(Fm)\Leftrightarrow U(Fm) = F\delta_{X}(U(Fm))\leq U'(Fm) = U'(Fm).$ Thus 
		$\mathcal{B'}\leq \mathcal{B}$.
	\end{proof}
	\begin{prop}\label{P1}
		Let $F : \mathcal{A}\longrightarrow \mathcal{C}$ be a faithful $\mathcal{M}$-fibration and $\mathcal{S}$ be a simple co-perfect syntopogenous structure on 
		$\mathcal{C}$ with respect to $\mathcal{M}$ i.e $\mathcal{S} = \{\sqsubset_{X}\}\in \bigwedge-INTORD$ . Then $c^{\sqsubset^{F}}(m) = \delta_{X}(c^{\sqsubset}(Fm))$ is an idempotent closure operator on $\mathcal{A}$
		with respect to $\mathcal{M}_{F}$. It is the largest closure operator on $\mathcal{A}$ for which $F$ is $(c^{\sqsubset^{F}}, c^{\sqsubset})$-continuous. 
	\end{prop}
	\begin{proof}
		It is easily seen that $c^{\sqsubset^{F}}$ is a closure operator for any simple co-perfect syntopogenous structure $\mathcal{S}$.  
		Now, $c^{\sqsubset^{F}}(c^{\sqsubset^{F}}(m))  =  c^{\sqsubset^{F}}(\delta_{X}(c^{\sqsubset}_{FX}(Fm)))
		= \delta_{X}(c^{\sqsubset}_{FX}(F\delta_{X}(c^{\sqsubset}(Fm))))
		= \delta_{X}(c^{\sqsubset}_{FX}(c^{\sqsubset}_{FX}(Fm)))
		= \delta_{X}(c^{\sqsubset}_{FX}(Fm)) = c^{\sqsubset^{F}}(m)$, thus 
		$c^{\sqsubset^{F}}$ is idempotent.
		$F$ is $(c^{\sqsubset^{F}}, c^{\sqsubset})$-continuous since, $\gamma_{X}(c^{\sqsubset^{F}}(m)) = c^{\sqsubset}(Fm)\Leftrightarrow Fc^{\sqsubset^{F}}(m) = c^{\sqsubset}(Fm)$.
		If $c'$ is another closure operator on $\mathcal{A}$ such that $F$ is $(c', c^{\sqsubset})$-continuous, then $Fc'_{X}(m)\leq c^{\sqsubset}(Fm)$. Thus 
		$c'_{X}(m) = \delta_{X}(F(c_{X}(m))\leq \delta_{X}(c^{\sqsubset}_{FX}(Fm)) = c^{\sqsubset^{F}}_{X}(m)$. 
	\end{proof}
	The closure operator in Proposition \ref{P1} was already obtained in  \cite {MR2122803} without use of the methods of syntopogenous structures. The interested
	reader will, in this book, find a number of examples for such closure.
	
	\thm
		Let $F\dashv G : \mathcal{C}\longrightarrow \mathcal{A}$ be adjoint functors and $\mathcal{B}$ be a base for a quasi-uniformity $\mathcal{U}\in $ QUnif$(\mathcal{C}, \mathcal{M})$. Assume that $G$ and $F$ preserve subobjects. Then 
		$\mathcal{B}^{\eta}_{X} = \{U^{\eta}\;|\;U\in \mathcal{B}_{FX}\}$ with $U^{\eta}(m) = \eta_{X}^{-1}(GU(Fm))$ for any $X\in \mathcal{A}$ is a base for a quasi-uniformity on $\mathcal{A}$. $\mathcal{B}^{\eta}$ is a base for the coarsest quasi-uniformity for which $F$ is $(\mathcal{U}^{\eta}, \mathcal{U})$-continuous.
	\endthm
	\begin{proof}
		Let us first note that $(U1)$, $(U2)$ and $(U4)$ are easily seen to be satisfied
		by adjointness.
		For $(U5)$, let $X\longrightarrow Y$ be a $\mathcal{A}$-morphism and $U^{\eta}\in \mathcal{U}^{\eta}$ for any $U\in \mathcal{U}_{Y}$. Then there is $V\in \mathcal{U}_{X}$ such that $f(V(m))\leq U(f(m))$.
		\begin{align*}
		\text{Thus}\;f(V^{\eta}(m)) &= f(\eta_{X}^{-1}(GV(Fm))))\\
		&\leq \eta_{Y}^{-1}(GFf)(GV(Fm))) &&\text{Lemma \ref{L22}} \\
		&\leq \eta^{-1}_{X}(G(Ff)(V(Fm)))\\
		&\leq \eta_{Y}^{-1}(GU((Ff)(Fm))) &&\text{$\mathcal{U}$-continuity of $Ff$} \\
		&= \eta_{Y}^{-1}(GU(Ff(m)))\\
		& = U^{\eta}(f(m)).\\
		\end{align*}
		$F$ is $(\mathcal{U}^{\eta}, \mathcal{U})$-continuous, since for any $U\in \mathcal{U}_{FX}$, $FU^{\eta}(m)\leq U(Fm)$ for any $X\in \mathcal{C}$.
		Let $\mathcal{B'}$ be a base for another quasi-uniformity $\mathcal{U}$ on $\mathcal{C}$ such that $F$ is $(\mathcal{U'}, \mathcal{U})$-continuous. Then for any $U^{\eta}\in \mathcal{B}_{X}^{\eta}$, there is $U'\in \mathcal{B'}_{X}$ such that $FU'(m)\leq U(Fm)$. 
		Thus $\eta_{X}(U'(m))\leq GFU'(m)\leq GU(Fm)\Rightarrow \eta_{X}(U'(m))\leq GU(Fm))\Leftrightarrow U'(m)\leq \eta_{X}^{-1}(GU(Fm)) = U^{\eta}(m)$, that is $\mathcal{U}^{\eta}\leq \mathcal{U'}$.
	\end{proof}
	If $\mathcal{A}$ is a reflective subcategory of $\mathcal{C}$, then $\mathcal{B}^{\mathcal{A}}$ and $\mathcal{B}^{\eta}$ are equivalent.
	\begin{prop}\label{P9}
		Let $F\dashv G : \mathcal{C}\longrightarrow \mathcal{A}$ be adjoint functors and $\mathcal{S}\in $ CSYnt$(\mathcal{C}, \mathcal{M})$. Assume that $G$ and $F$ preserves subobjects. 
		
		Then
		$\mathcal{S}^{\eta} = \{\sqsubset^{\eta}_{X}\;|\;\sqsubset_{FX}\in \mathcal{S}_{FX}\}$ 
		with $m\sqsubset^{\eta}_{X} n\Leftrightarrow \eta_{X}^{-1}(GU^{\sqsubset}(Fm))\leq n$ is a coperfect syntopogenous structure on $\mathcal{A}$. It is the coarsest syntopogenous structure for which $F$ is $(\mathcal{S}^{\eta}, \mathcal{S})$-continuous.
	\end{prop}
	\begin{prop}
		Under the assumptions of Proposition \ref{P9}, if $\mathcal{S}\in $ CSYnt$(\mathcal{C}, \mathcal{M})$ and simple i.e $\mathcal{S} = \{\sqsubset_{X}\}\in \bigwedge-$INTORD$(\mathcal{C}, \mathcal{M}) \cong ICL(\mathcal{C}, \mathcal{M})$. Then $c^{\sqsubset^{\eta}}_{X}(m) = \eta_{X}^{-1}(Gc_{FX}^{\sqsubset}(Fm))$ is an idempotent closure operator on $\mathcal{A}$. It is the largest closure operator for which $F$ is $(c^{\sqsubset^{\eta}}, c^{\sqsubset})$-continuous. 
	\end{prop}
	
	\section{Examples}
	\begin{enumerate}
		\item  Let $\bf QUnif_{o}$ be the category of $T_{o}$ quasi-uniform spaces and quasi-uniformly continuous maps with (surjective, embeddings)-factorization system. It is known 
		that $\bf bQUnif_{o}$ (see e.g \cite{MR1718995}), the category of bicomplete quasi-uniform spaces and quasi-uniformly continuous maps is an epi-reflective subcategory of $\bf QUnif_{o}$.
		Let $(F, \eta)$ be the bicompletion reflector into $\bf QUnif_{o}$.
		For any $(X, \mathcal{U})\in \bf QUnif_{o}$, $\eta_{X} : (X, \mathcal{U})\longrightarrow (\widetilde{X}, \widetilde{\mathcal{U}})$
		takes each $x\in X$ to its neighbourhood filter in the topology induced by the join of $\mathcal{U}$ and its inverse.
		It is known that $\eta_{X}$ is a quasi-uniform embedding. Details about this can be found in \cite{fletcher1982quasi}.
		Now, $\mathcal{B}^{F,\eta} = \{U^{F, \eta}\;|\;\widetilde{U}\in \widetilde{\mathcal{U}}_{\widetilde{X}}\}$ where 
		$U^{F, \eta} = \{(x, y)\in X\times X\;|\; (\eta_{X}(x), \eta_{X}(y))\in \widetilde{U}\}$ is a base for the quasi-uniform structure $\mathcal{U}^{F, \eta}$ on $X$.
		Since $\eta_{X}$ is  quasi-uniform embedding, $\mathcal{U}_{X}$ is the initial quasi-uniformity for which $\eta_{X}$ is quasi-uniformly continuous. 
		Thus  $\;\mathcal{U}^{F, \eta}_{X} = \mathcal{U}_{X}$.
		\item The category $\bf Unif$ of uniform spaces and quasi-uniformly continuous maps is coreflective  in $\bf QUnif$. Let $(G, \varepsilon)$ be 
		the coreflector into $\bf Unif$. For any $(X, \mathcal{U})\in \bf QUnif$, $\varepsilon_{X} : (X, \mathcal{U}\bigvee \mathcal{U}^{-1})\longrightarrow (X, \mathcal{U})$ is an identity map. Since $\mathcal{U}\bigvee \mathcal{U}^{-1}$ is the finest quasi-uniformity on $X$ for which $\varepsilon_{X}$ is quasi-uniformly continuous, $\mathcal{U}^{G, \varepsilon}_{X} = \mathcal{U}\bigvee \mathcal{U}^{-1}$  
		\item Consider $\bf TopGrp_{2}$ the category of Hausdorff topological groups and continuous group homomorphisms with the (surjective, injective)-factorization structure.
		We know from \cite{n1998general} that the category $\bf cTopGrp_{2}$ of complete Hausdorff topological groups (those topological groups which are complete with respect to the two-sided
		uniformity) is coreflective in $\bf TopGrp_{2}$. Let $(F,\eta)$ be the completion reflector into $\bf TopGrp_{2}$ and for any $(X, \cdot)\in \bf cTopGrp$, let $\beta(e)$ be the neighbourhood filter of the identity element $e$.
		For  all $U\in \beta(e)$, put $U_{c} = \{(x,\;y)\in X\times X\;: y\in xU\cap Ux\}$ so that 
		$\mathcal{B}^{c}_{X} = \{U^{c}\;|\;U\in \beta(e)\}$ is a base for the two-sided uniformity $\mathcal{U}^{c}$ on $(X; \cdot; \mathcal{T}\}$. Since $\eta_{X}$ is again an embedding of $(X,\cdot, \mathcal{T})\in \bf TopGrp_{2}$ into its completion $(\widetilde{X};\widetilde{\cdot}, \widetilde{\mathcal{T}})$, we have that  
		$\mathcal{U}^{F,\eta} = \mathcal{U}_{c}$.

		\item The forgetful functor $$F : \bf TopGrp\longrightarrow \bf Grp$$ is a mono-fibration. Thus by Proposition \ref{38}, every syntopogenous structure on $\bf Grp$ can be initially lifted to a syntopogenous structure on $\bf TopGrp$.
		\item  Consider the functors $G : \bf QUnif\longrightarrow \bf Top$ which sends every quasi-uniform space $(X, \mathcal{U})$ to the topological space $(X, G(\mathcal{U}))$ with 
		$G(\mathcal{U})$, the topology induced by $\mathcal{U}$, obtained by taking a base of neighbouhoods at a point $x$ the filter $\{U[x]\;|\;U\in \mathcal{U}\}$ where 
		$U[x] = \{y\in X\;:\;(x,y)\in U\}$ and $F : \bf Top\longrightarrow \bf Qunif$ which sends every topological space $(X, \mathcal{T})$ to the finest quasi-uniformity $\mathcal{U}$ on $X$
		with $G(\mathcal{U}) = \mathcal{T}$.
		It is known (see e.g \cite{MR1785843}) that $F$ is left adjoint to $G$. 
		For any $(X, \mathcal{T})\in \bf Top$, the unit $\eta_{X} : (X, \mathcal{T})\longrightarrow (X, GF(\mathcal{T}))$ is a continuous map where 
		$(X, GF(\mathcal{T}))$ is the set $X$ endowed with the topology of the finest quasi-uniformity $(X, F(\mathcal{T}))$.
		Now $\mathcal{S}_{(X, \mathcal{U})} = \{\sqsubset^{U}_{X}\;|\;U\in \mathcal{U}\}$  where $A\sqsubset^{U} B\Leftrightarrow U(A)\subseteq B$ for any $A, B\subseteq X$ is a co-perfect syntopogenous structure on $\bf Qunif$ for any $(X, \mathcal{U})\in \bf Qunif$.
		Let $(X, \mathcal{T})\in \bf Top$, $A\sqsubset_{X}^{\eta} B \Leftrightarrow \eta_{X}^{-1}(GU(FA))\subseteq B$ for any $U\in \mathcal{U}_{X}$.
		But $\eta_{X}^{-1}(GU(FA)) = \eta_{X}^{-1}(GU(A))$,  
		$\eta_{X}^{-1}(GU(A))$ is a neighbourhood of $A$ in $\mathcal{T}$. 
		Thus $\mathcal{S}_{X} = \{\sqsubset^{\eta}_{X}\;|\;X\in \bf Top\}$ with $A\sqsubset^{\eta}_{X} B \Leftrightarrow V\subseteq B$ where $V$ a is  neighbourhood of $A$ in $\mathcal{T}$
		so that $A\sqsubset^{\eta}_{X} B\Leftrightarrow A\subseteq O\subseteq B$ for some $O\in \mathcal{T}$.
		\item Let $\bf Top$ be the category of topological spaces and continuous maps with its (surjections, emdeddings)-factorization structure. It is well known that $\bf Top_{o}$, the category of
		$T_{o}$-topological spaces and continuous maps is a epi-reflective subcategory of $\bf Top$. Define $\mathcal{S}_{X} =\{\sqsubset_{Xo}\;|\;Xo\in \bf Top_{o}\}$ by $A\sqsubset_{X_{o}} B\Leftrightarrow \overline{A}\subseteq B$ for any $X_{o}\subseteq \bf Top_{o}$, $A, B\subseteq X_{o}$.
		Let $(F ,\eta)$ be the reflector into $\bf Top$. For any $X\in \bf Top$, $\eta_{X} : X\longrightarrow X/\sim$ takes each $x\in X$ to its equivalence class
		$[x] = \{y\in X\;|\;\overline{\{x\}} = \overline{\{y \}}\}$. Thus $\mathcal{S}_{X} = \{\sqsubset^{F, \eta}_{X}\;|\;X\in \bf Top\}$ with $A\sqsubset^{F, \eta}_{X} B\Leftrightarrow \eta_{X}^{-1}(\overline{\eta_{X}(A)})\subseteq B$ $A, B\subseteq X$.
		
\end{enumerate}

\bibliographystyle{abbrv}
\bibliography{references}
\addcontentsline{toc}{chapter}{References}
\end{document}